\documentclass{amsart}

\usepackage{amsmath}
\usepackage{amsfonts}
\usepackage{amssymb}
\usepackage{amsthm}

\newcommand{\Z}{{\bf{Z}}}
\newcommand{\Q}{{\bf{Q}}}

\newcommand{\C}{{\bf{C}}}
\newcommand{\F}{{\bf{F}}}
\newcommand{\Fp}{{\F}_p}
\newcommand{\OO}{\mathcal{O}}
\newcommand{\HH}{\mathcal{H}}
\newcommand{\Ell}{{\rm{Ell}}}

\newcommand{\Mid}{\mid\!\!}
\newcommand{\miD}{\!\!\mid}

\newcommand{\ra}{\rightarrow}

\newcommand{\nn}{\newline \noindent}
\newcommand{\noi}{\noindent}
\newcommand{\bigo}{{\rm O}}
\newcommand{\rem}{{\rm rem}}

\newcommand{\comment}[1]{}

\theoremstyle{plain}
\newtheorem{lem}{Lemma}[section]

\newtheorem{prop}[lem]{Proposition}

\newtheorem{stat}[lem]{Statement}

\begin{document}

\title[Constructing elliptic curves]{Constructing elliptic curves with a known number 
of points over a prime field}

\author{Amod~Agashe}
\address{Department of Mathematics, University of Texas, Austin, Texas, USA}
\email{amod@math.utexas.edu}

\author{Kristin~Lauter}
\address{Microsoft Research,
         One Microsoft Way,
         Redmond, WA 98052, USA.}
\email{klauter@microsoft.com}

\author{Ramarathnam~Venkatesan}
\address{Microsoft Research,
         One Microsoft Way,
         Redmond, WA 98052, USA.}
\email{venkie@microsoft.com}

\thanks{We would like to thank the anonymous
referees of an earlier version of the paper for pointing out mistakes
and making suggestions. The first author would like to thank
D.~Kohel, H.~Lenstra, F.~Morain and J.~Vaaler for several 
useful conversations, and M.S.R.I., I.H.E.S., and M.P.I. for
their generous hospitality.  The second author would like to thank
R.~Schoof for pointing out the paper \cite{cnst}, and P.~Montgomery, D.~Bernstein,
A.~Stein, and H.~Cohn for useful discussions. 
}

\date{}


\renewcommand{\subjclassname}{\textup{2000} Mathematics Subject Classification}


\begin{abstract}
Elliptic curves with a known number of points over a given 
prime field~$\F_n$ are often needed for use in cryptography.
In the context of primality proving, Atkin and Morain suggested the use 
of the theory of complex multiplication to construct such curves. One of 
the steps in this method is the calculation of a root modulo~$n$ of
the Hilbert class polynomial $H_D(X)$ for a fundamental discriminant $D$. 
The usual way is to compute 
$H_D(X)$ over the integers and then to find the root modulo $n$.
We present a modified version of the Chinese remainder theorem (CRT)
to compute $H_D(X)$ modulo~$n$ directly
from the knowledge of $H_D(X)$ modulo enough small primes. 
Our complexity analysis suggests that asymptotically our algorithm is an 
improvement over previously known methods.
\end{abstract}

\maketitle

\section{Introduction} \label{section:intro}

In order to use elliptic curves in cryptography, one often needs to construct
elliptic curves with a known number of points over a given 
prime field. One way of doing this is to randomly pick elliptic curves and
then to count the number of points on the curve over the prime field, 
repeating this until the desired number of points is found.
Atkin and Morain \cite{atmor} pointed out that instead, one can use the theory of complex 
multiplication to construct elliptic curves with a known number of points.
Although at present it may still be more efficient to count
points on random curves, we hope that improving the complex multiplication
method will eventually yield a more efficient algorithm.  In some 
situations, using complex multiplication methods is the only practical
possibility (e.g. if the prime is too large for point-counting to be
efficient yet the discriminant of the imaginary quadratic field is
relatively small).
This paper provides a new version of the complex multiplication method.

Suppose $n$ is an integer, usually a prime or a pseudo-prime,
and one wants to construct an elliptic curve modulo~$n$ along with the 
number of points on that curve modulo~$n$. 
One of the steps in the complex multiplication method 
is the calculation of the Hilbert class polynomial 
$H_D(X)$ modulo~$n$ for a certain fundamental discriminant $D$. 
The usual way to do this is to compute 
$H_D(X)$ over the integers and then to reduce it modulo~$n$.
Atkin and Morain proposed computing $H_D(X)$ as an integral polynomial 
by listing all the relevant binary quadratic forms, associating to each form an
algebraic integer, evaluating the
$j$-function at each of those as a floating point integer with sufficient precision, 
and then taking the product and rounding the coefficients to nearest
integers. Let $d = \Mid D \miD$.
If we use the estimate given by formula~(\ref{upperbound}),
then in view of~\cite[\S5.10]{ll}, the
computation of~$H_D(X)$ by this method takes time $\bigo(d^2 (\log d)^2)$.

In~\cite[\S4]{cnst}, the authors suggested computing $H_D(X) \bmod p$ 
for sufficiently many small primes $p$ and then using
the Chinese remainder theorem (CRT) to compute $H_D(X)$ as a polynomial with 
integer coefficients.  
In this paper we use a modified version of CRT 
to compute $H_D(X)$ modulo~$n$ directly (knowing
$H_D(X) \bmod p$ for sufficiently many small primes $p$),
{\it without} computing its coefficients as integers.
We also give the mathematical justification and details of
the (usual) CRT method, which were omitted
in~\cite[\S4]{cnst} and also correct their erroneous complexity analysis.
By avoiding the computation of the coefficients of $H_D(X)$ as integers,
we obtain an algorithm 
with asymptotically shorter running time as $d$ gets large.
Also, both CRT approaches require less precision of computation than the Atkin-Morain approach. 

Our complexity analysis in Section~\ref{section:algo} shows that,
when $d$ is large, with high probability, 
the running time of one of the versions of our algorithm is
$$\bigo(d^{3/2} (\log d)^{10} + d (\log d)^2 \log n + \sqrt{d} (\log n)^2 ),$$
which is better than the Atkin-Morain method
when $d$ is sufficiently large (roughly speaking, bigger than
$(\log n)^2$). Our algorithm has a step in common
with the (usual) CRT method, which takes time $\bigo(d^{3/2} (\log d)^{10})$,
and for the other step, our algorithm takes time
$$\bigo(d (\log d)^4  + d (\log d)^2 \log n + \sqrt{d} (\log n)^2 ),$$ 
while the (usual) CRT method takes time 
$$\bigo( d (\log d)^2 \log n + d^{3/2} (\log d)^4  ).$$
Thus we obtain an improvement over the (usual) CRT  method when $d$ is 
greater than $(\log n)^2$.

Note that in~\cite{atmor}, the authors suggest that using Weber polynomials
works better in practice than using Hilbert polynomials. At the moment, we do not have
a generalization of our algorithm which works with Weber polynomials.
The use of Weber polynomials only reduces the number of digits by a constant,
hence will only change the time taken by a constant factor   
independent of~$d$, (see~\cite[p.409]{cohen}), so the asymptotic complexity 
estimates remain the same.

Note also that we only focus on one step of the complex multiplication
algorithm, the computation of the Hilbert class polynomial,
The other time-consuming step is the 
computation of a root of~$H_D(X)$ modulo~$n$, which (by~\cite[\S5.10]{ll})
takes time $\bigo(d (\log n)^3)$. The relative size of~$d$ and~$n$ will determine
which of these two steps will dominate (when we use our algorithm to
compute $H_D(X)$ modulo~$n$). 

It is not clear how our method compares to existing
methods computationally. While we did some examples
(reported in Section~\ref{examples}), they involved small discriminants,
where existing methods are already very fast. The purpose of
this paper is to suggest a new version of the complex multiplication
method and to present a complexity analysis, leaving the task of
efficient implementation for the future.

The paper is organized as follows:
in Section~\ref{section:cm_method}, we give a brief description of the 
complex multiplication method for generating elliptic curves.
In Section~\ref{section:algo}, we give an outline of our
algorithm and discuss its complexity. 
In Sections~\ref{section:hdxmodp} and~\ref{section:couv}, 
we explain the details of some of the steps of the algorithm. 
Finally in Section~\ref{examples}, we give some
examples of our method.


\section{Complex multiplication method} \label{section:cm_method}

We briefly review the complex multiplication method, referring the reader to~\cite{atmor}
and~\cite{silv2} for details.
Suppose we are given a prime $n$, and a non-negative number $N$ in
the Hasse-Weil interval $[n+1-2\sqrt{n},n+1+2\sqrt{n}]$. 
We want to produce an elliptic curve $E$ over $\F_n$ with $N$ 
points over $\F_n$:
$\#E(\F_n) = N =n+1-t,$
where $t$ is the trace of the Frobenius endomorphism of $E$ over $\F_n$.
We set $$D=t^2-4n.$$
The Frobenius endomorphism of $E$ has characteristic polynomial
$x^2-tx+p$, and its roots lie in $\Q(\sqrt{D})$.  It is standard to associate the Frobenius 
endomorphism with a root of this polynomial.  If $t\ne 0$, then $E$ is not supersingular,
in which case $R$, the endomorphism ring of $E$, is an order in the ring
of integers of $K=\Q(\sqrt{D})$ (\cite[Thm 3.1.b]{silv1}).  
For simplicity of the algorithm, we will want to assume
that $R$ is $\OO_K$, the full ring of integers in $K$.  Recall that 
a negative integer~$D$ is said to be a {\em fundamental discriminant} if it is not divisible
by any square of an odd prime and satisfies $D \equiv 1 \bmod 4$
or $D \equiv 8, 12 \bmod 16$.  If $D$ is a fundamental discriminant,
then $R$ is automatically equal to $\OO_K$, since then the Frobenius endomorphism generates the
full ring of integers and is contained in the endomorphism ring.
Our results can be generalized to orders in the ring of integers, but the 
algorithm will become more complicated. We will assume throughout this paper that $D$ is a 
fundamental discriminant.  In particular, this means that the simplest version of our
algorithm only works for those choices of $n$ and $N$ such that this condition on $D$ is met.

The {\em Hilbert class polynomial}~$H_D(X)$ is defined as:
\begin{eqnarray} \label{hdx}
H_D(X) = \prod \bigg(X - j \bigg(\frac{-b+\sqrt{D}}{2a} \bigg) \bigg),
\end{eqnarray} 
where the product ranges over the set of $(a,b) \in \Z \times \Z$ such that
$ax^2 + bxy + cy^2$ is a primitive, reduced, positive definite binary
quadratic form of discriminant~$D$ for some $c \in \Z$,
and $j$ denotes the modular invariant. The degree of $H_D(X)$ is equal
to $h$, the class number of $\OO_K$.  
It is known that $H_D(X)$ has integer coefficients.
The equivalence between isomorphism classes of elliptic curves over $\bar{\Q}$
with endomorphism ring equal to $\OO_K$ and primitive, reduced, positive 
definite binary quadratic forms of discriminant~$D$ allows us to interpret
a root of this polynomial as the j-invariant of an elliptic curve having
this endomorphism ring.  Since our goal is to find such an elliptic
curve modulo~$n$, it suffices to find a root $j$ of ~$H_D(X)$ modulo~$n$. 

Assuming $j \neq 0, 1728$, the required elliptic curve is recovered as 
the curve with Weierstrass equation (assume $n \ne 2,3$)
$$y^2=x^3+3kx+2k,$$
where $$k=\frac{j}{1728-j}.$$
The number of points on the elliptic curve is either
$n+1+t$ or~$n+1-t$, and one can easily check which one
it is by raising randomly chosen points to one of the possible group orders.

\section{Our algorithm and its complexity} \label{section:algo}

\subsection{Overview of the algorithm} \label{algo-overview}

As before, let $D$ be a fundamental discriminant and let $d = \Mid D \miD$.
Let $K = \Q(\sqrt{D})$ , let $\OO_K$ denote the ring
of integers of~$K$, and let $h$ denote the class number of $\OO_K$.      
Let $B$ be an upper bound on the size of the coefficients of~$H_D(X)$
given by the formula in Section \ref{estimates}.
Let $n$ be a given prime number.

Here is our algorithm for computing $H_D(X) \bmod n$; it comes
in two versions, Version~A and Version~B, which differ
only in Step~(1) below:

%
%

\vspace{0.1in}
\noi {\bf Step~(0)} Compute $h$ and~$B$.
Compute $h$ using any of the standard algorithms
(e.g., see~\cite[\S5.4]{cohen})
and compute~$B$ using formula~(\ref{Bformula}) in Section~\ref{estimates}.
Fix a small real number $\epsilon > 0$ (e.g. $\epsilon = 0.001$),
and let $M = B / (1/2-\epsilon)$.
\vspace{0.1in}

\noi {\bf Step~(1)} Compute $H_D(X)$ modulo sufficiently many small primes:

\vspace{0.1in}

\noi {\bf Version~A:} This can be used whenever $d \not \equiv 7 \bmod 8$. 

\vspace{0.1in}

\noi (a)
Generate a collection of distinct primes $p$, each satisfying
$4p = t^2 + d,$  for some integer $t$.
Generate enough primes $p$
so that the product of all the primes exceeds the bound $B$
(or slightly exceeds $2B$, see the remark after Example 6.1).


\noi (b) 
For each~$p$ in~$S$, consider a set of representatives for the
$\overline{\F}_p$-isomorphism classes of elliptic curves over~$\F_p$, and
count the number of $\F_p$-points on each representative.  
In practice, we take as a representative the model $$y^2=x^3+3kx+2k,$$
where $k=\frac{j}{1728-j},$ and $j$ runs through all possible values in
$\F_p$ (except $0$ and $1728$, which can be handled separately if necessary).
We then form the set $S_p$ consisting of all the $j$-invariants
such that the corresponding curve has $p+1+t$ or $p+1-t$ points.    
There are exactly $h$ such $j$ values, by
Prop.~\ref{prop:hdxmodp} and Prop.~\ref{prop:points} below
(or by \cite[p. 319]{cox}).  Alternatively, for each representative, 
we could pick random points~$P$ on~$E$
and check if $(p+1)P = t P$ (or $(p+1)P = - t P$). This would rapidly
filter out almost all of the candidates, and point-counting could be
used to check the remaining ones. 

\noi (c) For each prime $p$ in~$S$, we form the polynomial $H_D(X) \bmod{p}$
by multiplying together the factors $(X-j)$, where $j$ is in the set $S_p$. 
This is also justified by Prop.~\ref{prop:hdxmodp} and
Prop.~\ref{prop:points} below.

\vspace{0.1in}

\noi {\bf Version~B:} This can be used for any~$d$; however, we expect
it to be more difficult to implement.

\vspace{0.1in}

\noi Version~B is exactly like Version~A except that we allow slightly
more general primes when forming the set~$S$ in Step~(a).
We allow all primes $p$ such that
$4p = t^2 + u^2 d,$ for some integers~$t$ and~$u$. 
We again generate enough primes $p$
so that their product exceeds the bound $B$;
call the resulting set of primes~$T$.
Then for each $p$ in~$T$, we compute the endomorphism ring for each 
$\overline{\F}_p$ isomorphism class of elliptic curves over $\F_p$
using the algorithm in~\cite{kohel} (we use the same representatives
for the isomorphism classes as in Version~A, Step~(b) above).
We then form the set $T_p$ consisting of all the 
$j$-invariants such that the corresponding curve has endomorphism ring
isomorphic to~$\OO_K$.
The class number of $\OO_K$ is $h$, so there are exactly $h$ such $j$ values.

Finally, as in Version~A, Step~(c), for each prime $p$ in~$T$, 
we form the polynomial $H_D(X) \bmod{p}$
by multiplying together the factors $(X-j)$, where $j$ is in the set $T_p$. 

\vspace{0.1in}
\noi {\bf Remark 3.1.}
Note that in Version~B, when we allow more general primes~$p$ 
such that $4p = t^2+u^2d,$ where $u > 1$, it is not sufficient to use 
point-counting to find the desired collection of elliptic curves. 
In that case, point-counting would produce the set of all elliptic curves with 
endomorphism ring equal to an order in $\OO_K$ containing the order of index $u_i$.
In this paper, we assumed that $d$ was square-free, but
to generalize our algorithm to non-square-free $d$, it would be 
necessary to work with Version~B of the algorithm.  The number and size 
of the primes required to implement the two versions does not seem to be much 
different in practice (see the remark after Example 6.2).  The main advantage 
to Version~A is that it is easy to implement because there are many 
point-counting packages available.  
The main advantage to Version B is that it will generalize to work for all~$d$.

\vspace{0.1in}
\noi {\bf Step~(2)} Lift to $H_D(X) \bmod n$: 

\vspace{0.1in}
\noi Use the modified chinese remainder algorithm of Section~\ref{section:couv} to compute
each coefficient of $H_D (X) \bmod n$ using the values of the coefficients of
$H_D(X) \bmod p$ computed in Step~(1).  This step can be parallelized.

\subsection{Complexity anaylsis} \label{complexity-analysis}

In our complexity analysis, we assume that if $a$ and~$b$ are two integers,
then their addition takes time $\bigo(\log a + \log b)$, their
multiplication takes time $\bigo(\log a \log b)$, and the division
of the greater by the smaller takes time $\bigo(\log a \log b)$.
This can certainly be achieved by current algorithms; in fact,
one can do better, but we will stick to our model of computation
for the sake of simplicity and comparison (the complexity estimate for the 
Atkin-Morain algorithm, $\bigo(d^2)$, 
given in~\cite{ll} does not assume fast arithmetic either). 
The steps mentioned below are numbered as in Section~\ref{algo-overview}.

\vspace{0.1in}

\noi {\bf Step~(0)} 
 According to~\cite[\S5.4]{cohen}, the computation of~$h$
can be done in time $\bigo(d^{1/4})$, or in time $\bigo(d^{1/5})$ 
assuming the generalized Riemann Hypothesis, and $B$ is computed from the formula
given in Section~\ref{estimates}.

\vspace{0.1in}

\noi {\bf Step~(1)} We do the analysis only for Version~A.

\noi (a)
By the discussion in \S~\ref{estimates},
with high probability, 
the size of~$S$ is $\bigo(\frac{\log B}{\log d})$, and
each $p \in S$ is $\bigo((\log B)^2)$; for the purposes
of the complexity analysis, we will assume this happens
(this makes our complexity analysis ``probabilistic'').

\noi (b) 
The best implementations of elliptic curve point-counting algorithms 
currently run in time $\bigo ((\log p)^5)$ (\cite{schoof}), perhaps assuming
fast arithmetic, although this will not affect the power of $d$ in our overall
complexity estimate.  This step is repeated $p$ 
times, so this step will take time $\bigo (p(\log p)^5)$.
Finally, since the step is repeated for every prime in~$S$,
the total time taken will be $\bigo((\log B)^3 (\log \log B)^5/\log d )$.
In Section \ref{estimates}, we estimate $\log B$ in terms of $d$ as
$ \log(B) = \bigo (\sqrt{d} (\log d)^2)$.  Using this estimate,
the time taken for this step in terms of $d$ is
$\bigo(d^{3/2} (\log d)^{10})$, ignoring $\log \log d$ factors.
We should be able to speed up this step in practice by using the alternative
suggested above to avoid counting points on each curve modulo~$p$.

\noi (c) 
The number of terms in the product used to compute
$H_D(X) \bmod{p}$ is~$h$ and each coefficient
is between zero and~$p$, so this can be done in time 
$\bigo (h^2 (\log p)^2)$, i.e., $\bigo(d (\log d)^2)$.
Since the step has to be repeated for every~$p\in S$,
the total time taken is $\bigo(d^{3/2} (\log d)^3)$.

\vspace{0.1in}
Overall, the total time taken by Step~$1$ in this version
is $\bigo(d^{3/2} (\log d)^{10})$.

\vspace{0.1in}
\noi {\bf Step~(2)} 
As will be explained in Section~\ref{section:couv},
the time taken by the modified chinese remainder algorithm to compute
all the coefficients of $H_D (X) \bmod n$ is
$$\bigo(d (\log d)^4  + d (\log d)^2 \log n + \sqrt{d} (\log n)^2 ).$$

\vspace{0.1in}

Our algorithm differs from the one in~\cite[\S4]{cnst} mainly in Step~(2).
As shown in Section~\ref{section:couv}, 
if one uses the ordinary Chinese remainder theorem
to find $H_D(X)$ and then reduces modulo~$n$, as proposed in~\cite[\S4]{cnst},
then the complexity of this procedure would be
$$\bigo( d (\log d)^2 \log n + d^{3/2} (\log d)^4  ),$$
which is not as good as our method in Step~(2) 
when $d$ is large (roughly speaking, bigger than $(\log n)^2$).

On the other hand, for primality proving as in~\cite{atmor},
one wants a small discriminant; in fact, in 
\cite[\S5.10]{ll} they assume $d = \bigo( (\log n)^2)$. 
In that case, it is clear that our algorithm is an improvement over the one in~\cite[\S4]{cnst}
only if $\log B$ is bigger than~$\log n$, i.e., if 
the coefficients of $H_D(X)$ are large compared to~$n$. 

The overall complexity of our algorithm, assuming 
Statement~\ref{estimate-stat}, is
$$\bigo(d^{3/2} (\log d)^{10} + d (\log d)^2 \log n + \sqrt{d} (\log n)^2 ).$$

\subsection{Some estimates needed for the complexity analysis}
\label{estimates}

We need an estimate for the size of $B$, i.e., an 
upper bound for the size of the coefficients of the class polynomial.
As is explained in \cite[p.~42]{atmor}, we may take 
\begin{eqnarray} \label{Bformula}
B = \left(
  \begin{array}{c}
  h \\
  \lfloor h/2 \rfloor 
  \end{array} \right) 
{\rm exp}\bigg({\pi \sqrt{d}} \sum \frac{1}{a} \bigg),
\end{eqnarray}
where the sum in the above expression is taken over the set of integers~$a$
such that $ax^2 + bxy +cy^2$ is a primitive, reduced, positive definite 
binary quadratic form of discriminant~$D$ for some integers~$b$ and~$c$
(the set of $a$'s is finite).  
This bound comes from the product of all the roots times the largest 
binomial coefficient.  

Note that by the corollary in \cite[Chap XVI, \S4]{langant},
we have $\log h \sim \log(\sqrt{d})$ as $d \ra \infty$
(recall that the regulator of a quadratic imaginary field is one).
This means that for any positive real number~$\epsilon'$, we have
$d^{1/2-\epsilon'} \leq h \leq d^{1/2+\epsilon'}$
when $d$ is big enough. For the sake of simplicity 
in our analysis, we will assume $h \sim \sqrt{d}$.

We will soon need a lower bound on the size of~$\log B$.
By \cite[Lem.~5.3.4(1)]{cohen}, $a \leq \sqrt{d/3}$.
Thus $\sum \frac{1}{a} \geq h \sqrt{\frac{3}{d}}$, and
the latter is asymptotically a constant bigger than~$1$.
Thus there is a constant~$c>1$ such that $\log B$ 
is greater than $c \sqrt{d}$
for $d$ large enough.

To get an upper bound for $\log B$ in terms of $d$, we
estimate  $\sum \frac{1}{a}$ 
using the argument in~\cite[p.~711]{ll}. They observe that
there cannot be too many $a$'s  that are ``small'', since the 
the number of reduced forms $(a,b)$ with a fixed $a$ is bounded 
by $\tau(a)$, the number of positive divisors of $a$. 
So certainly an overestimate for the sum $\sum \frac{1}{a}$ 
is given by $\sum_{a=1}^d \frac{\tau(a)}{a}$.  This in turn can
be written as a telescoping sum plus an error term:
$$\sum_{a=1}^d \frac{\tau(a)}{a} = 
\sum_{a=1}^d (\sum_{u=1}^a \tau(u))(\frac{1}{a}-\frac{1}{a+1}) +
\frac{1}{d+1}\sum_{a=1}^d \tau(a).$$
The sum $\sum_{a=1}^d \tau(a)$ can be estimated as $d\log d$ plus some lower order terms
(see \cite[Thm 8.28, p. 393]{nzm}).  
So the first term can be estimated via the integral
$$\int_{a=1}^d \frac{\log a}{a} da = \frac{(\log d)^2}{2},$$
and the second term is less than $\log d$.
This observation leads to the estimate $$ \sum \frac{1}{a} \le \bigo ((\log d)^2),$$
(see also~\cite[p.~324]{crandall-pomerance}). 
In fact, much better estimates for $\sum \frac{1}{a}$ should be possible,
and it looks like a better bound is being assumed in
the complexity analysis for the Atkin-Morain algorithm given
by~\cite{ll}, since they seem to assume that $\log(B) = \bigo (\sqrt{d})$, but we will
stick with our estimate for our analysis.

Since the middle binomial coefficient is clearly less than the sum of all of the
binomial coefficients, which is $2^h$, 
we see that $$B \le 2^h e^{\pi\sqrt{d}(\log d)^2}.$$
So throughout the paper, we use the estimate
\begin{eqnarray} \label{upperbound}
\log(B) = \bigo (\sqrt{d} (\log d)^2) = \bigo(h (\log h)^2).
\end{eqnarray}

An important consideration for accurately assessing the running time of 
our algorithm is the relative size of the small primes found in Step~(1). 
Consider the following statement:
\begin{stat} \label{estimate-stat}
If $d \not \equiv 7 \bmod 8$, then the procedure of finding primes 
in Version~A of Step (1) terminates, and  
the size of the set $S$ is $\bigo(\frac{\log B}{\log d})$ and
each $p \in S$ is $\bigo((\log B)^2)$.

\end{stat}
We expect that the statement above is true with high probability
when~$d$ is large enough. The main idea for Statement 3.1 was
suggested to us by an anonymous referee. We now give a heuristic
argument to support our expectation, some of the details
of which were explained to us by J.~Vaaler.

By the prime number theorem, the probability that a randomly chosen
positive integer~$m$ is prime is $1/(\log m)$.  For a given $d$,
and randomly chosen $t$, we want to say that a number of the form $(t^2+d)/4$
looks like a randomly chosen integer, so that we can claim that the probability that it is
prime is $1/\log((t^2+d)/4)$. 

If $d \equiv 3 \mod 8$, say $d = 8k + 3$, and  
if $t$ is odd, say $t = 2\ell + 1$, then $(t^2 + d)/4 = 
\ell (\ell+1) + 2k +1$ is an odd integer.
If $d \equiv 4 \bmod 16$, say $d = 16 k + 4$,
and if $t$ is a multiple of~$4$, 
say $t = 4 \ell$, then $(t^2 + d)/4 = 4 \ell^2 + 4k +1$
will be an odd integer.
If $d \equiv 8 \bmod 16$ (the only possibility left), 
say $d = 16 k + 8$,
and if $t$ is even, say $t = 2 \ell$, then $(t^2 + d)/4 = \ell^2 + 4k +2$
will be an odd integer provided $\ell$ is odd.
So for any~$d$, for a random choice of an integer~$t$,
with probability at least~$1/4$, the rational number $(t^2 + d)/4$ 
will be an odd integer
(i.e., $(t^2 + d)/4$ will be an integer that need not 
necessarily be composite).
So we will assume that the probability that it is prime 
(provided it is an odd integer) is indeed $1/\log((t^2+d)/4)$. 

Now let $c_1$ and~$c_2$ be two positive integers such that $c_1 < c_2$.
Let $S_1$ denote the set
$$S_1 = \{(t^2+d)/4 : t \in \Z, c_1 \log B \leq t \leq c_2 \log B, 
(t^2+d)/4 \mbox{{\rm\ is prime}} \}.$$
The size of the set 
$\{(t^2+d)/4 : t \in \Z, c_1 \log B \leq t \leq c_2 \log B\}$ is $(c_2 - c_1) \log B$, 
and roughly one-fourth of the elements of this set 
are integers. Moreover, among those which are integers,
we are assuming that the probability that an element $(t^2+d)/4$ is
prime is $1/\log((t^2+d)/4)$. Thus with high probability,
the following statement is true for large~$d$:
\newline $(*)$ The size of the set~$S_1$ is between
$ \frac{1}{2} \left \lfloor \frac{(c_2 - c_1) \log B}{4 \log(c_2 \log B)} 
\right \rfloor$ 
and  $ 2 \left \lfloor \frac{(c_2 - c_1) \log B}{4 \log(d/4)} \right \rfloor$. 
\newline 
We will assume that $(*)$ is indeed true for the rest of this section
(so everything below holds only with high probability).

If $p \in S_1$, then 
$$p \geq \frac{(c_1 \log B)^2 + d}{4} > \frac{ (c_1 \log B)^2}{4}.$$
Thus $$\sum_{p \in S_1} \log p > [2 (\log \log B) + \log (c_1^2/4)]
\frac{(c_2 - c_1) \log B}{4 \log(c_2 \log B)}.$$
By choosing $c_1$ and $c_2$ appropriately (say $c_2=12$ and $c_1=4$), we see
that when $d$ is large enough (so that $\log B > c_2$),
$\sum_{p\in S_1} \log p > \log B$
and hence $\prod_{p\in S_1} p > B$.

Now let $S_2$ denote the set
$$S_2 = \{(t^2+d)/4 : t \in \Z, 0 \leq t \leq c_2 \log B, 
(t^2+4)/d \mbox{{\rm\ is prime}} \}.$$
Putting $c_1=0$ in statement~$(*)$, we see that the size of~$S_2$ 
will be $\bigo(\frac{\log B}{\log d})$. 
\newline
Also, $\prod_{p\in S_1} p > B$, since the set~$S_2$
contains the set~$S_1$. Furthermore, if $p \in S_2$, then 
$$p < ((c_2 \log B)^2 + d)/4.$$ 
Since $d$ is $\bigo((\log B)^2)$,
we see that $p$ is $\bigo((\log B)^2)$.
Finally (assuming statement $(*)$ holds), the set~$S$ can be 
chosen to be a subset of
the set~$S_2$; from this, Statement~\ref{estimate-stat} follows.

\section{Computing $H_D(X) \bmod p$ for small primes~$p$}
\label{section:hdxmodp}

In this section, we prove that Step~1~ of our algorithm is a valid way
to compute $H_D(X) \bmod p$. 
The same strategy for this step was used in~\cite[\S4]{cnst},
but it was not justified there, and the distinction between 
Versions~A and ~B was blurred.

As in the introduction, let $D$ be a fundamental discriminant 
and let $H_D(X)$ denote the Hilbert class polynomial.
Let~$H$ denote the Hilbert class field of~$K=\Q(\sqrt{D})$.
and let $p$ be a rational prime that splits completely in~$H$,
i.e., splits into principal ideals in~$K$, which means that
$4p = t^2 - D u^2$ for some integers~$u$ and~$t$. 

Let $\Ell(D)$ denote the set of isomorphism classes 
of elliptic curves over~$\C$ with complex multiplication  
by~$\OO_K$ (i.e., whose ring of endomorphisms over~$\C$ is
isomorphic to~$\OO_K$). Then an equivalent way of defining
the Hilbert class polynomial is as follows: 
\begin{eqnarray} \label{hdx2}
H_D(X) = \prod_{[E] \in \Ell(D)} (X - j(E)),
\end{eqnarray}
where, if $E$ is an elliptic curve, then $j(E)$ denotes its $j$-invariant. 

Let $\Ell'(D)$ denote the set of~$\overline{\F}_p$-isomorphism classes 
of elliptic curves over~$\F_p$ with 
endomorphism ring (over~$\overline{\F}_p$) isomorphic to~$\OO_K$.

\begin{prop} \label{prop:hdxmodp}
With notation as above,
\begin{eqnarray} \label{hdxmodp}
H_D(X) \bmod p = \prod_{[E'] \in \Ell'(D)} (X - j({E'})).
\end{eqnarray}
\end{prop}

\begin{proof}
Let $\beta$ be a prime
ideal of the ring of integers of~$H$ lying over~$p$.
It follows from the discussion in the proof of Thm.~14.18 
on p.~319-320 of~\cite{cox} that in each class~$i$ in
$\Ell(D)$, we can write down an elliptic curve~$E_i$ such that
$E_i$ is defined over~$H$ and $E_i$ has good reduction modulo~$\beta$
(in fact, \cite{cox} gives a collection of such elliptic curves,
denoted~$E_c$; we just pick one such~$E_c$ for each class); 
denote the reduction modulo~$\beta$ of $E_i$ by 
$\widetilde{E_i}$. 
Since $p$ splits completely in~$H$,
$\widetilde{E_i}$ is defined over~$\Fp$, 
as opposed to an extension of~$\Fp$.
Also, by~\cite[Chap~13, Thm.~12(ii)]{lang}
(or \cite[Thm.~14.16]{cox}), each $\widetilde{E_i}$ has 
endomorphism ring (over~$\overline{\F}_p$) isomorphic to~$\OO$.
This gives us a map~$\phi$ from $\Ell(D)$ to $\Ell'(D)$.
Since we assume that $p$ splits in $K$, then
by \cite[Thm.~13.21]{cox}, if two elliptic curves 
have distinct $j$-invariants, then the reductions
modulo~$\beta$ of these $j$-invariants are distinct,
i.e., the map~$\phi$ is injective.
By the Deuring lifting theorem \cite[Chap.~13, Thm.~14]{lang}
(or \cite[Thm.~14.16]{cox}) 
this map is also a surjection.

From the definition of~$j(E)$ in terms of the coefficients of the Weierstrass
equation of~$E$, it is easy to see that 
$$H_D(X) \bmod p = \prod_{[E_i] \in \Ell(D)} (X - j(\widetilde{E_i})).$$
Hence, from the discussion above,
$$H_D(X) \bmod p = \prod_{[E'] \in \Ell'(D)} (X - j({E'})).$$

\comment{
Let $S$ denote the set of elliptic curves over~$\C$
with complex multiplication by~$\OO$,
and let $S'$ denote the set of elliptic curves over~$\F_p$
with complex multiplication by~$\OO$.
First, by~\cite[Chap~13, Thm.~12(ii)]{lang}
(or \cite[Thm.~14.16]{cox}), 
reduction modulo~$\beta$ gives a map from $S$ to~$S'$.
Since we assume that $p$ splits in $K$, then
by \cite[Thm.~13.21]{cox}, if two elliptic curves 
have distinct $j$-invariants, then the reductions
modulo~$\beta$ of these $j$-invariants are distinct,
i.e., we have an induced injective map from
$\Ell(D)$ to $\Ell'(D)$.
By the Deuring lifting theorem \cite[Chap.~13, Thm.~14]{lang}
(or \cite[Thm.~14.16]{cox}) 
this map is also a surjection.
}

\end{proof} 
 
\begin{prop} \label{prop:points} Recall that $D$ is a fundamental discriminant.
Suppose $p$ is a prime and $x \ne 0$ is an integer  such that \mbox{$4p = x^2 -D$}.
Let $E'$ be an elliptic curve over~$\F_p$.
Then $[E'] \in \Ell'(D)$ if and only if
$\# E'(\F_p)$ is either $p+1-x$ or $p+1+x$.
\end{prop}

\begin{proof}
Suppose $\# E'(\F_p)$ is either $p+1-x$ or $p+1+x$.
Let $t$ denote the trace of the Frobenius
endomorphism of~$E'$.  Then $t=x$ or $t=-x$.
In either case, the discriminant of the characteristic
polynomial of the Frobenius endomorphism
is $t^2 - 4p = x^2 - 4p = D$. Let ${\rm End}(E')$ denote
the endomorphism ring of~$E'$. 
Since $D$ is square-free, the subring $R$ of ${\rm End}(E')$ generated by
the Frobenius endomorphism is $\OO$, and at the 
same time ${\rm End}(E')$ is contained in the ring
of integers of the quotient field of~$R$.
Hence ${\rm End}(E') = \OO$, i.e., $[E'] \in \Ell'(D)$. 

Conversely, suppose $[E'] \in \Ell'(D)$,
and let $t$ denote the trace of the Frobenius
endomorphism of~$E'$. Suppose the Frobenius
endomorphism generates a subring of index~$u$
in ${\rm End}(E')$, the endomorphism ring of~$E'$. 
Then the characteristic polynomial of the Frobenius endomorphism
has discriminant $u^2 D$, hence $4p = t^2 - u^2 D$.
But we know $4p = x^2 - D$, so by \cite[Ex.~14.17]{cox},
$t = x$ or $t = -x$. Hence 
$\# E'(\F_p)$ is either $p+1-x$ or $p+1+x$.
\end{proof}

\section{A modification of the Chinese remainder theorem} \label{section:couv}

\subsection{The algorithm and its complexity} \label{couveignes-algo}

This section follows~\cite[\S2.1]{couv} closely,
which in turn is based on~\cite[\S4]{mont-silv}; the only addition
is a more detailed complexity analysis.

The problem we consider is as follows: for some positive integer~$\ell$ 
we are given a collection
of pairwise coprime positive integers~$m_i$ for $i=1,2, \ldots, \ell$.
For each~$i$, we are also given an integer~$x_i$ with $0 \leq x_i < m_i$. 
In addition, we are given 
a small positive real number~$\epsilon$. Finally, we are told that there is 
an integer~$x$ 
such that $|x| < (1/2 -\epsilon) \prod_i m_i$
and $x \equiv x_i \bmod m_i$ for each~$i$; clearly such an integer~$x$
is unique if it exists. The question is to compute $x \bmod n$,
for a given positive integer~$n$.

Define
\begin{eqnarray}
\label{eqnM} M &=& \prod_i m_i \\
\label{eqnMi} M_i &=& \prod_{j \neq i} m_j = M/m_i \\
\label{eqna} a_i &=& 1/M_i \bmod m_i, \ \ \ \ \ 0\leq a_i < m_i.
\end{eqnarray}
Then the number $z = \sum_i a_i M_i x_i$ is congruent to~$x$
modulo~$M$. 
Hence, if $r = \big \lfloor \frac{z}{M} + \frac{1}{2} \big \rfloor$, 
then $x = z - r M$. 
So $x \bmod n = z \bmod n - (r \bmod n)(M \bmod n)$;
the point is that we can calculate $r \bmod n$ without calculating~$z$,
as we now explain. 
 From the fact that
$x = z - rM$ and $|x| < (1/2 -\epsilon) M$, it follows
that $\frac{z}{M} + \frac{1}{2}$ is not within~$\epsilon$ of
an integer. Hence, to calculate~$r$, one only has find an approximation~$t$
to $z/M$ such that $|t - z/M| < \epsilon$, 
and then round~$t$ to the nearest integer. Such an approximation~$t$
can be obtained from 
\begin{eqnarray} \label{eqnzM}
\frac{z}{M} = \sum_i \frac{a_i x_i}{m_i},
\end{eqnarray} 
where the calculations are done using floating point numbers.

If $a$ and~$b$ are two integers, then 
let $\rem(a,b)$ denote the remainder of the Euclidean division
of~$a$ by~$b$; we will assume that it takes time $\bigo(\log a \log b)$
to calculate $\rem(a,b)$ and ${\rm gcd}(a,b)$. 


From the discussion above, we obtain the following algorithm:
\nn (i) Compute~$a_i$'s, for each~$i$, using~(\ref{eqna}): this takes
time $\bigo(\sum_i ( \sum_j (\log m_j \log m_i) + \ell (\log m_i)^2
+ (\log m_i)^2)) = \bigo((\log M)^2 + \ell \sum (\log m_i)^2)$.
\nn (ii) Compute $\rem(M,n)$ using~(\ref{eqnM}): this will take
time $\bigo(\sum_i (\log m_i \log n) + \ell (\log n)^2 )
=\bigo(\log n \log M + \ell (\log n)^2 )$.
\nn (iii) Compute $\rem(M_i,n)$ for each~$i$ by
dividing $\rem(M,n)$ by~$m_i$ modulo~$n$: this will take time
$\bigo(\ell (\log n)^2)$ (in our application, $m_i$ will be much
lesser than~$n$), and can be parallelized.
\nn (iv) Compute~$r$: 
In~(\ref{eqnzM}), every term in the sum has to be calculated
to precision $\epsilon/\ell$, hence the calculation of each
term takes time $\bigo( (\log(\ell/ \epsilon))^2 )$. 
In the application to computing $H_D(X) \bmod n$, 
we can take $\epsilon$ to be an arbitrary small number 
and taking $M = B/(1/2-\epsilon)$.
Then the calculation of all the terms in~(\ref{eqnzM}) 
will take total time $\bigo(\ell (\log \ell)^2)$ 
and the addition in~(\ref{eqnzM}) of $\ell$ numbers with 
precision~$\epsilon/\ell$ will take time $\bigo(\ell \log \ell)$.
\nn (v) 
Output $\rem(x,n) = $
\begin{eqnarray} \label{eqnfinal}
\rem \Bigg( \bigg( \rem \Big( \sum_i (\rem(a_i \cdot x_i,n) \cdot \rem(M_i,n) ), n \Big)
- \rem(r,n) \cdot \rem(M,n) \bigg), n \Bigg).
\end{eqnarray}
The various substeps in step (v) and the time taken for each are as follows:
\nn (a) Calculation of $\rem(a_i \cdot x_i,n)$ and $\rem(M_i,n)$
for all~$i$:
takes time $\bigo(\sum_i ((\log m_i)^2 + (\log m_i)(\log n)))$.
\nn (b) Computing the product of $\rem(a_i \cdot x_i,n)$ and
$\rem(M_i,n)$ for all~$i$: takes time $\bigo(\ell (\log n)^2)$.
\nn (c) Performing the sum in~(\ref{eqnfinal}) and taking remainder
modulo~$n$: this involves about $\ell$ additions of
integers of size up to~$\ell n^2$, which takes time 
$\bigo(\ell \log(\ell n^2))$ and taking the remainder takes time 
$\bigo((\log n)(\log(\ell n^2))$.
\nn (d) Calculation of $\rem(r,n) \cdot \rem(M,n)$: The size of $r$ is
about $\sum m_i$, hence this substep takes time 
$\bigo((\log n) \log (\sum m_i) + (\log n)^2)$.
\nn (e) Subtraction operation and taking remainder: takes time 
$\bigo(\log n)$ and $\bigo((\log n)^2)$ respectively.

In Section~\ref{section:algo}, we use this algorithm to lift 
$H_D(X) \bmod p$ for
$p \in S$ to $H_D(X) \bmod n$ one coefficient at a time.
Note that steps (i), (ii), and (iii) above are common to the lifting
of all the coefficients,
and only step (iv) and (v) have to be repeated for each coefficient.

In the notation of Section~\ref{section:algo}, 
the $m_i$'s are the elements of~$S$, and so, assuming
Statement~\ref{estimate-stat},
we see that $m_i$'s are $\bigo( (\log B)^2 )$ and $\ell$ 
is~$\bigo( \log B/\log\log B )$.
Using this, and the estimates from~\S~\ref{estimates},
we see that the most time consuming steps
are Step~(i), which takes time $\bigo(d (\log d)^4)$,
and Steps~(v-a) and~(v-d) repeated $h$~times, which take
time $\bigo(d (\log d)^2 \log n)$ and $\bigo(\sqrt{d} (\log n)^2 )$
respectively.
\comment{
time taken (asymptotically) for calculations 
in steps (v)(a), (v)(b) and (v)(c), performed $h$  
times, (once for each coefficient of $H_D(X)$) will dominate all other steps.
We first calculate the overall complexity of these 3 steps using the estimates
$\ell \approx \log B$ and $p_i \approx \log B$.  This is true for the 
set of primes with no conditions of the type imposed in Step~(1)(a) of
 our algorithm.
In that case the complexity would be:
$$\bigo( ( \log p \log B \log \log B +(\log p + \log B)(\log p+\log \log B))),$$
and this needs to be repeated $h$ times.
This expression can be simplified using the fact that $\log\log B < \log p$, 
yielding
$$\bigo( \log p (\log B) (\log \log B) + (\log p)^2),$$
which should be performed $h$ times.
}


\vskip .1 truein

\subsection{Complexity of the usual Chinese Remainder Algorithm}

If we are to use the naive Chinese remainder theorem for the
problem stated at the beginning of Section~\ref{couveignes-algo}, 
then we calculate 
\begin{eqnarray} \label{eqncrt}
z 
= \rem \bigg( \Big(\sum_i a_i \cdot x_i \cdot M_i \Big), M \bigg), 
\end{eqnarray}
and then reduce $z$ modulo~$n$.

The steps involved are as follows:
\nn (i) Compute~$a_i$'s, for each~$i$, using~(\ref{eqna}): 
this takes
time $\bigo(\sum_i ( \sum_j (\log m_j \log m_i) + \ell (\log m_i)^2
+ (\log m_i)^2)) = \bigo((\log M)^2 + \ell \sum (\log m_i)^2)$.
\nn (ii) Calculation of $a_i \cdot  x_i \cdot M_i$ for all~$i$: takes time 
$\bigo ( \sum_i (\log m_i)(\log M))$. 
\nn (iii) Performing the sum in~(\ref{eqncrt}): this involves $\ell$ additions
of integers of size up to~$\ell m_i^2 M$, 
hence takes time $\bigo(\ell \log(\ell m_i^2 M))$.
\nn (iv) Calculating the outer ``\mbox{$\rem$}'' in~(\ref{eqncrt}):
takes time $\bigo( (\log M) \log (\ell m_i^2 M))$.
\nn (v) Reducing $z$ modulo~$n$: takes time $\bigo((\log M)(\log n))$.

In the context of lifting $H_D(X) \bmod p$ to $H_D(X)$ and then reducing
$H_D(X)$ modulo~$n$, only steps (ii) -- (vi) have to be repeated for
each coefficient. 
In the notation of Section~\ref{section:algo}, 
the $m_i$'s are the elements of~$S$, and so, assuming
Statement~\ref{estimate-stat},
we again have that $m_i$'s are $\bigo( (\log B)^2 )$ and $\ell$ 
is~$\bigo( \log B/\log\log B )$.
Using this, and the estimates from~\S~\ref{estimates},
we see that the most time consuming steps are Steps~(ii) and~(iv), 
each of which take total time $\bigo( d^{3/2} (\log d)^4)$ 
and Step~(v), which takes total
time $\bigo( d (\log d)^2 \log n)$.

\comment{
Again in the notation of Section~\ref{section:algo}, 
$m_i=p_i$.  For the sake of comparison, again assume that $\ell$ and~$p_i$
are $\bigo( \log M)$, and $M = B$. Thus the time taken for the
calculation of step~(ii) above,
performed $h$ 
times (once for each coefficient of~$H_D(X)$), 
will dominate to give a complexity of 
$\bigo( (\log B)^2 (\log \log B) )$ to be performed $h$ times.
}

From this analysis, we see that our modified Chinese remainder 
algorithm will be asymptotically more efficient than the usual one 
when $\log B > \log n $, which will certainly be the case in our context 
whenever $d \ge (\log n)^2$ (certainly, when $n > B$, the modified
version is no better than the usual Chinese remainder algorithm).

\section{Examples} \label{examples}

In this section we present several examples to illustrate 
our algorithm.  Throughout these examples, we used the software package PARI, which is
available at 
\begin{center}
{\tt http://www.parigp-home.de}
\end{center}

\subsection{$D= -59$}
\subsubsection{Atkin-Morain Method}
Since here we are dealing with a very small discriminant, we can easily
compute the minimal polynomial over the integers directly by finding all the
reduced, positive definite, primitive, binary quadratic forms with 
discriminant $-59$ and then evaluating $j(\tau)$ for the corresponding $\tau$ 
with sufficiently high precision.  The class number of $\Q(\sqrt{-59})$ is three,
and the three binary quadratic forms are 
$$(a,b,c) = (3,1,5), (3,-1,5), (1,1,15).$$
The corresponding algebraic integer is 
$$\tau_{(a,b,c)} = \frac{-b + \sqrt{b^2-4ac}}{2a}.$$
We expect the absolute value of the largest of the $j(\tau)$
to be roughly $e^{\pi\sqrt{59}} \approx e^{24}$.
Evaluating the product $$(x-j(\tau_1))(x-j(\tau_2))(x-j(\tau_3))$$
with enough significant digits and rounding the coefficients
to integers, we find the class polynomial:
$$H_D(x) =  x^3 + 30197678080x^2 - 140811576541184x + 374643194001883136.$$
Here $28$ decimal digits of precision are required using the 
package pari ($19$ digits of precision are not enough).

\subsubsection{Chinese Remainder type algorithms}
To implement our algorithm for this example, we set the bound
$B$ equal to $e^{41}$ to be bigger than the largest 
coefficient of $H_D(x)$.  This estimate comes from the product
of the three $j$ values, whose absolute value we expect to be
roughly $$e^{\pi\sqrt{59}(1+\frac{1}{3}+\frac{1}{3})}.$$

We find the following list of $7$ small primes which 
are of the form $(t^2 -D)/4$ for some integer~$t$:
$$17,71,197,521,827,1907,3797,5417$$
and whose product exceeds~$B$.  For each prime~$p$ in the list,
we loop through the $p-1$ possible $j$-values.  For each possible
$j$-value, we count the number of points on a curve over $\F_{p}$ 
with that $j$-value using a version of Schoof's algorithm
(we use a version available on the web by Mike Scott: 
{\tt ftp://ftp.compapp.dcu.ie/pub/crypto/sea.cpp}).
If the curve has either $p+1+t$ or $p+1-t$ points,
with $t^2 = 4p-59$, then we keep that $j$-value in a list $S_p$.
At the end of the loop, we will have $h$ $j$-values in the list $S_p$,
where $h$ is the degree of $H_D(x)$.  Then the polynomial $H_D(x) \mod {p}$
is formed as the product over $j \in S_p$ of $(x-j)$.

Here is a table summarizing the results for this example:

\vskip .1 truein
\begin{center}


\begin{tabular}{|c|c|c|c|}

\hline
$p$ & $t$ & $j \in S$ & $H_D(x) \mod {p}$  \\

\hline
$17$ & $3$ & $j=2,7,13$  &     $x^3+12x^2+12x+5$ \\

\hline
$71$ & $15$ & $j=51,54,67$ &   $x^3+41x^2+62x+11$ \\

\hline
$197$ & $27$ & $j=71,195,130$   & $ x^3+195x^2+160x+139$ \\

\hline
$521$ & $45$ & $j=103,366,367$  & $ x^3+206x^2+379x+510$ \\

\hline
$827$ & $57$ & $j=97,498,554$ &    $x^3+505x^2+824x+196$ \\

\hline
$1907$ &$87$ & $j=24,915,1613$ & $x^3+1262x^2+1432x+1045$ \\

\hline
$3797$ & $123$ & $j=70,958,2381$ & $x^3+388x^2+1114x+1584$ \\

\hline

\end{tabular}
\end{center}

\vskip .1 truein
\noi
{\bf Usual Chinese Remainder routine}
\vskip .1 truein

Here is a short routine in the algebraic number theory package PARI
to compute the polynomial $H_D(x)$ with integer coefficients using the usual
Chinese Remainder Theorem.  It takes as input the coefficients of
$H_D(x)$ modulo the small primes $p$.

\vskip .1 truein

l=7;                     (number of small primes)

h=degree;  (degree of the hilbert class polynomial)

m=[17,71,197,521,827,1907,3797];  (list of small primes)

M=prod(i=1,l,m[i]);    (M=17*71*197*521*827*1907*3797)

log(M);

invm = vector(l,i,M/m[i]);

a=vector(l,i,Mod(1/invm[i],m[i]));

modcoeff = $[[12,41,195,206,505,1262,388],[12,62,160,379,824,1432,1114],$

\noindent
$[5,11,139,510,196,1045,1584]]$;  (list of coefficients modulo small primes)

z=vector(h,j,Mod(sum(i=1,l,lift(a[i])*invm[i]*modcoeff[j][i]),M));

\vskip .1 truein
\noi
{\bf Modified Chinese Remainder routine}

\vskip .1 truein

For our algorithm, we input in addition the prime $n$ such that we want
to determine $H_D(x) \mod {n}$. Here is a short routine in PARI
to compute the polynomial $H_D(x)$ with coefficients modulo $n$ using 
our modified version of the Chinese Remainder Theorem.

\vskip .1 truein

n=prime;   (the prime where we want the curve in the end)



r=vector(h,j,round(sum(i=1,l,(lift(a[i])*modcoeff[j][i]/m[i]))))

finalcoeff=vector(h,j,sum(i=1,l,

$\quad$ $\quad$ lift(a[i])*modcoeff[j][i]*Mod(invm[i],n))-Mod(r[j],n)*Mod(M,n))

\vskip .1 truein

Note that the precision required for this computation is almost trivial
(the minimum value to set the precision in PARI is $9$ significant digits).

\vskip .1 truein
\noi
{\bf n=141767}
\vskip .1 truein

Here is an example where we use our algorithm to find the class polynomial 
modulo $n$.  Note that $4n = 753^2 - D$, so we will construct a curve over $\F_n$
with $142521 = n+1 + 753$ points.  The output of our Modified Chinese Remainder
routine is:
$$[{\rm Mod}(31177, 141767), {\rm Mod}(73152, 141767), {\rm Mod}(48400, 141767)].$$
Note that this corresponds to the class polynomial that we found using the Atkin-Morain method 
reduced modulo $n$:
$$X^3+31177X^2+73152X+48400.$$
Taking the root $j=118481 \bmod n$, we get the elliptic curve 
$$y^2=x^3+39103x+120580.$$ It has $142521$ points as desired.

\vskip .1 truein
\noi 
{\bf Remark 6.1.}
Actually, the third coefficient in this example had to be re-computed 
because there was a rounding problem.  The constant term of the class 
polynomial over the integers is slightly more than half the product 
of the small primes.  The problem in this example can be solved in a clean 
way by adding one more prime to the algorithm, since in fact our algorithm requires
the product of the small primes to slightly exceed $2B$ 
by an amount depending on the choice of epsilon: $B/(1/2-\epsilon)$.


\subsection{$D=-832603$}

\vskip .1 truein

\noindent
The Algorithms and Parameters for Secure Electronic Signatures document put out
by the EESSI-SG (European Electronic Signature Standardisation Initiative 
Steering Group) recommends using elliptic curves with class number of the
endomorphism ring at least equal to $200$.  Here is an example
with class number equal to $96$.  Let $$n=100959557.$$ 
Note that $4n=20075^2-D$, and so we will construct a curve over $\F_n$ with
$N = 100979633 = n+1 +20075 $ points.  

We have that $D$ is square-free and $\Q(\sqrt{D})$ has
class number $h= 96$, which is small compared to the square root of $|D|$,
$$\sqrt{|D|} \approx 912.$$
According to the estimates, the largest coefficient of the class polynomial is
bounded by $e^{5368}$.  This comes from the fact that 
$$\sum_{i=1}^{96}\frac{1}{a} \approx 1.85$$ and the middle binomial coefficient 
is roughly $e^{64}$.

\subsubsection{Atkin-Morain method}

We can obtain the class polynomial using the algorithm of Atkin and Morain 
with 3000 digits of precision (2332 digits should suffice): 

\tiny
$x^{96} +$ 
8986950689916460612768050899826095370126160959774067006607495722714787327536405

11195940426329493962250363608918506814518954357512108376324309765509813261009595

72615396095460780845684222178125741276615369508754546593823796954336290438786044

54401141760139087536761069637319825798963352735300683356996983745448647386647867

26063390570575243929901907691553896874023573783820406248132762600812249711372047

55882861009465622598491768081847062179144325418471334571707958715942702715493314

88237482404374709398037956562818329691426448791666838198286258706039015068098946

24510449977402232596376577346850486922319170621800447828468015555155662062177391

08385797919357857853408831828384120143178961174846657318648760182117564137653818

31734687436523641791869559811287475164662560558340565130954332294988968060971888

71815593515878469206432064483048317524773444570122581660831541350800516869161291

01483247617224314616156733489349276043330450686852025326165481636562782630791850

43524347061886083145402858558832786452505054211954992588893518489408407045712834

80364209087452918765509915544167886763955395126621398677472529126929317764001654

07674073078383580568650075515962375620983618886988248866522341997936320370535130

5474956970365974518712304022211825509601280000$x^{95} \dots$ 

\normalsize

The coefficients of this polynomial are so big that it would take about $30$
pages just to write the polynomial down in that font size.  

\subsubsection{Modified Chinese Remainder Algorithm}
To obtain the class polynomial using our method, we first make a
list of primes $p$ which are each of the form $4p = t^2 -D$ for some integer~$t$. 

\vskip .1 truein

\noindent
{\bf List of small primes}

\vskip .1 truein

\noindent
\tiny
[208207, 208223, 208261, 208283, 208333, 208391, 208457, 208493, 208657, 
208907, 208963, 209021, 210131, 210407, 210601, 210803, 210907, 211231, 211457, 
211573, 211691, 211811, 211933, 212057, 212183, 212573, 212843, 212981, 213263, 
213407, 213553, 214003, 214631, 215123, 215461, 215983, 216523, 217081, 217271, 
217463, 218453, 218657, 219071, 219281, 219707, 220141, 220361, 220807, 221261, 
221723, 221957, 222193, 222913, 223403, 224153, 224921, 226241, 226511, 226783, 
227611, 228457, 229321, 229613, 230203, 230501, 232643, 233591, 233911, 235211, 
235541, 236207, 236881, 237563, 240371, 241093, 241823, 245593, 245981, 246371, 
247553, 248351, 248753, 249563, 249971, 250793, 251623, 253307, 253733, 254161, 
255023, 255457, 256771, 257657, 258551, 259001, 259453, 259907, 260363, 261281, 
263611, 264083, 265511, 266957, 268913, 273943, 274457, 274973, 275491, 276011, 
278111, 279173, 279707, 280243, 281321, 282407, 283501, 284051, 286831, 287393, 
290233, 292541, 293123, 294293, 298451, 299053, 303323, 304561, 307691, 308323, 
310231, 315407, 320041, 321383, 322057, 324773, 326143, 326831, 328213, 329603, 
333821, 335957, 339557, 340283, 343943, 344681, 348401, 349903, 350657, 351413, 
352931, 356761, 364571, 366953, 367751, 368551, 369353, 373393, 375841, 379133, 
382457, 387503, 388351, 397811, 398683, 399557, 401311, 403957, 404843, 405731, 
411101, 414721, 415631, 416543, 417457, 418373, 419291, 421133, 422057, 424841, 
426707, 429521, 436157, 437113, 441923, 444833, 445807, 448741, 450707, 453671, 
456653, 457651, 458651, 460657, 466723, 468761, 474923, 475957, 480113, 481157, 
482203, 483251, 484301, 486407, 487463, 490643, 491707, 495983, 499211, 503543, 
504631, 507907, 510101, 511201, 513407, 514513, 515621, 520073, 523433, 527941, 
531343, 538201, 539351, 541657, 543971, 545131, 548623, 550961, 555661, 556841, 
560393, 568751, 571157, 573571, 579641, 582083, 593171, 598151, 606943, 608207, 
610741, 612011, 615833, 620957, 622243, 623531, 626113, 637831, 639143, 640457, 
644411, 647057, 648383, 652373, 653707, 655043, 659063, 665803, 668513, 672593, 
678061, 679433, 682183, 687707, 689093, 691871, 696053, 707293, 712961, 718661, 
720091, 724393, 727271, 728713, 730157, 731603, 735953, 738863, 740321, 741781, 
744707, 759457, 763921, 766907, 771401, 772903, 777421, 783473, 788033, 789557, 
794141, 797207, 800281, 806453, 814213, 817331, 823591, 825161, 829883, 831461, 
834623, 839381, 844157, 845753, 853763, 861823, 865061, 866683, 874823, 883013, 
897881, 901207, 902873, 906211, 911233, 912911, 914591, 916273, 921331, 924713, 
934907, 940031, 950333, 952057, 955511, 957241, 958973, 967663, 969407, 971153, 
972901, 974651, 976403, 978157, 983431, 997583, 1006493, 1015453, 1020853, 
1028081, 1031707, 1035341, 1040807, 1042633, 1048123, 1055471, 1066553, 1068407, 
1075843, 1077707, 1081441, 1083311, 1088933, 1094573, 1107803, 1115407, 1119221, 
1124957, 1130711, 1132633, 1134557, 1136483, 1138411, 1140341, 1146143, 1151963, 
1161703, 1167571, 1187261, 1195193, 1197181, 1201163, 1209151, 1217171, 1221193, 
1223207, 1225223, 1227241, 1231283, 1241423, 1278341, 1286633, 1288711, 1290791, 
1294957, 1301221, 1309601, 1311701, 1315907, 1318013, 1328573, 1330691, 1334933,
1337057, 1347707, 1351981, 1362701, 1377793, 1379957, 1382123, 1388633, 1392983,
1403893, 1406081, 1410463, 1417051, 1419251, 1425863, 1430281, 1432493, 1434707]

\vskip .1 truein

\normalsize
\noindent
The list contains $410$ primes. Their product is roughly $e^{5379}$, which
exceeds the bound $2B$ as desired.

To illustrate the algorithm, we find the class polynomial
modulo the largest prime on the list $p=1434707$. Note that
$4p = 2215^2 -D$. By counting 
the number of points on a representative for each isomorphism class of 
elliptic curves over $\F_p$, 
we found the following list of $96$ ~$j$-values such that the associated 
elliptic curve has $p+1 \pm 2215$ points over $\F_p$.

\vskip .1 truein

\noindent
{\bf $j$-values for $p=1434707$:}

\vskip .1 truein

\tiny
\noindent
[28534, 29664, 39989, 50559, 58497, 61669, 87155, 97333, 120663, 153566, 158121, 164378, 
182440, 199741, 210115, 218108, 219599, 
237389, 257474, 289215, 317239, 333891, 335757,
365925, 381504, 395862, 403801, 449952, 482780, 485134, 487074, 511916, 527120, 
543027, 574978, 583669, 584091, 585813, 595906, 642664, 644346, 653188, 654512, 655573, 696063,
698345, 699985, 702445, 705943, 710770, 721309, 738498, 759603, 780978, 795085, 816076,
821241, 869331, 871700, 889175, 897281, 902226, 923156, 924382, 980018, 1022428,
1033432, 1057121, 1079631, 1093031, 1101285, 1129437, 1154957, 1161878, 1175298,
1185913, 1186864, 1199076, 1205398, 1231078, 1252451, 1279055, 1281872, 1286184,
1312922, 1327236, 1334297, 1352254, 1352769, 1364919, 1368722, 1381024, 1410659,
1426507, 1428519, 1431597]

\vskip .1 truein

\small
\noindent
We find that $$H_D(X) \bmod{p} = $$
$X^{96}+1163995X^{95}+922656X^{94}+700837X^{93}+1079920X^{92}+466732X^{91}+154378X^{90}+399013X^{89}+744868X^{88}+1140439X^{87}+238431X^{86}+439229X^{85}+1168335X^{84}+1088371X^{83}+1065323X^{82}+923089X^{81}+370237X^{80}+418673X^{79}+26462X^{78}+1186790X^{77}+577727X^{76}+1026750X^{75}+1311499X^{74}+42221X^{73}+1226509X^{72}+1302356X^{71}+1205738X^{70}+706055X^{69}+916474X^{68}+870490X^{67}+940463X^{66}+779702X^{65}+543453X^{64}+1023692X^{63}+985646X^{62}+734246X^{61}+744646X^{60}+754597X^{59}+67621X^{58}+394070X^{57}+801259X^{56}+1203063X^{55}+1415480X^{54}+182257X^{53}+358715X^{52}+659376X^{51}+343711X^{50}+472997X^{49}+545620X^{48}+578548X^{47} +223638X^{46}+281011X^{45}+170375X^{44}+514817X^{43}+327182X^{42}+506290X^{41}+550176X^{40} +157534X^{39}+1257296X^{38}+1245604X^{37}+311058X^{36}+532467X^{35}+601208X^{34}+1069781X^{33}+52757X^{32}+508590X^{31}+247205X^{30}+1293507X^{29}+1089763X^{28}+326605X^{27}+46947X^{26}+1147567X^{25}+884035X^{24}+535907X^{23}+1164336X^{22}+952400X^{21}+1245681X^{20}+348341X^{19}+43230X^{18}+1201679X^{17}+486702X^{16}+360056X^{15}+28756X^{14}+1068784X^{13}+993753X^{12}+790102X^{11}+436946X^{10}+37636X^9+459204X^8+1185717X^7+644728X^6+1031301X^5+384651X^4+380850X^3+1358865X^2+1127134X+401105 \bmod{p}.$

\normalsize
\vskip .1 truein
\noindent
The class polynomials modulo the other $409$ primes are not included here.
This polynomial indeed corresponds to the reduction modulo $p$ of the class polynomial
found using the Atkin-Morain algorithm.

\vskip .1 truein
\noi
{\bf Remark 6.2.}
In this example, we find that allowing primes~$p$ such that  
$4p = u^2 + v^2d$ (i.e., using Version~B) does not help much.
The size of $v$ is constrained by the desire to keep the primes small; here, 
$v$ must satisfy $v \le 2$ to avoid getting larger primes.  Allowing 
$v=2$, we still need a list of length $410$ primes to exceed the bound.
The black art of balancing the size of the primes with the number of primes 
required is not the subject of this paper, but at least in this example Version~B seems
no better in this regard.

\end{document}